\documentstyle[12pt]{article}
\newcommand{\const}{\mathop{\rm const}\limits}

\newcommand{\mod}{\mathop{\rm mod}\limits}

\newcommand{\argmax}{\mathop{\rm argmax}\limits}

\newcommand{\argmin}{\mathop{\rm argmin}\limits}

\newcommand{\grad}{\mathop{\rm grad}\limits}

\begin{document}

\begin{center}

{\bf Non \ - \ asymptotic exponential bounds for }\\

\vspace{3mm}

 {\bf MLE deviation under minimal conditions}\\

\vspace{3mm}

  {\bf via classical and generic chaining methods.}\\

\vspace{4mm}

{\sc By Ostrovsky E., Rogover E.}\\

\vspace{3mm}

{\it Department of Mathematics and Statistics, Bar \ - \ Ilan University,
59200, Ramat Gan, Israel.}\\
e \ - \ mail: \ galo@list.ru \\
e \ - \ mail: \ eugeny@soniclynx.com\\

\vspace{3mm}

{\it Department of Mathematics and Statistics, Bar \ - \ Ilan University,
59200, Ramat Gan, Israel.}\\
e \ - \ mail: \ rogovee@gmail.com \\

\vspace{4mm}

   {\sc ABSTRACT}\\

\end{center}

\vspace{2mm}

 {\it In this paper non \ - \ asymptotic exact exponential estimates
 are derived (under minimal conditions) for the tail of deviation of the
 MLE distribution in the so \ - \ called natural terms: natural function,
 natural distance,
 metric entropy, Banach spaces of random variables,
 contrast function, majorizing measures or, equally, generic chaining. } \\

\vspace{3mm}


\vspace{3mm}


\vspace{3mm}

{\bf 1. Introduction. Notations. Statement of problem.}\\

\vspace{3mm}

 Let $ (\Omega,\cal{M},{\bf P} ) $ be a probability space with the
expectation $ {\bf E}, \ \Omega = \{ \omega \},
 \ (X, \cal{A}, \mu) $ be a measurable space with sigma \ - \ finite non \ - \ trivial measure $ \mu, \ (\Theta, \tau) = \Theta = \{\theta \} $ be arbitrary separable local compact topological space equipped by the ordinary Borelian sigma \ - \ field,
 $ \cal{F} $
  $ = \{ f \}, \ f = f(x,\theta) \ $ be a {\it family } of a {\it strictly positive} $ ( \mod \ \mu) $ probabilistic densities:

$$
\forall \theta \in \Theta \ \Rightarrow \int_X f(x,\theta) \ d \mu = 1,
$$

$$
\mu \{ \cup_{\theta \in \Theta} \{x: \ f(x,\theta) \le 0 \} \} = 0,
$$
continuous relative to the argument $ \theta $ for almost all values $ x;
x \in X. $ \par

{\it We premise also the following natural condition of the identifying: }

$$
\forall \theta_1, \theta_2 \in \Theta, \theta_1 \ne \theta_2 \ \Rightarrow \ 
\mu \{x: \ f(x, \theta_1) \ne f(x, \theta_2) \} > 0.
$$

 Let further $ \theta_0 \in \Theta $ be some fixed value of the parameter
$ \theta. $ We assume that $ \xi = \xi(\omega) $ is a random variable (r.v)
(or more generally random vector) taking the values in the space $ X $ with
the density of distribution $ f(x, \theta_0) $ relative the measure $ \mu: $

$$
{\bf P} (\xi \in A) = \int_A f(x,\theta_0) \ d \mu, \ A \in {\cal M}. \eqno(1.0)
$$

 The statistical sense: the r.v. $ \xi $ is the (statistical) observation
(or observations) with density $ f(x,\theta_0), $
 where the value $ \theta_0 $
is the "true", but in general case unknown value of the parameter $ \theta. $ \par
 We denote as usually by $ \hat{\theta} $ the Maximum Likelihood Estimation
(MLE) of the parameter $ \theta_0 $ based on the observation $ \xi: $

$$
\hat{\theta} = \argmax_{\theta \in \Theta} f(\xi, \theta),
$$
or equally

$$
\hat{\theta} = \argmax_{\theta \in \Theta} L(\xi, \theta) =
\argmin_{\theta \in \Theta}(- L(\xi, \theta)) \eqno(1.1)
$$
where the function

$$
 L = L(\xi, \theta) \stackrel{def}{=} \log \left[ f(\xi, \theta)/f(\xi, \theta_0) \right]
$$
is called the {\it contrast function}, in contradiction to the function
 $ \theta \to f(\xi,\theta) $ or $ \theta \to \log f(\xi, \theta), $
which is called ordinary {\it Likelihood } function. \par
 Denote also

 $$
 L^{(0)} = L^{(0)}(\xi, \theta) = L(\xi, \theta) \ -
\ {\bf E} L(\xi, \theta) = L(\xi, \theta) \ - \ {\bf E}_{0} \ L(\xi, \theta).
 $$

 In the case if $ \hat{\theta} $ is not unique, we understand as
$ \hat{\theta} $ any but measurable value $ \hat{\theta} $
satisfying the condition (1.1).\par

 Let now $ r = r(\theta) = r(\theta, \theta_0), \ \theta \in \Theta $ be some
(measurable) numerical non \ - \ negative {\it risk}, or {\it deviation}
 function, i.e. such that

$$
r(\theta, \theta_0) \ge 0, \
r(\theta, \theta_0) = 0 \ \Leftrightarrow \theta = \theta_0,
$$
not necessary to be distance, i.e. it can not satisfy the triangle inequality.
 We denote for arbitrary positive value $ v $ the probability of the deviation
in the $ r(\cdot,\cdot) $ sense of $ \hat{\theta} $ from the true value
$ \theta_0: $

$$
W(v) \stackrel{def}{=} {\bf P}( \ r(\hat{\theta}, \theta_0) > v), \eqno(1.2)
$$
which is needed for the construction of confidence region for the unknown
parameter $ \theta_0 $ in the $ r(\cdot, \cdot) $ sense. \par

\vspace{3mm}

{\bf Our goal of this paper is non \ - \ asymptotical estimation of the
function $ W = W(v) $ under minimal and natural conditions for sufficiently greatest values } $ v; \ v \ge v_0 = \const > 0, \ (v >>1). $ \par

 Offered here estimations are some generalizations of the paper \cite{Golubev1}.
See also \cite{Bagdasarova1}, \cite{Birge1}, \cite{Ibragimov1}, \cite{Jensen1}, \cite{Geer1} and reference therein. \par

 The paper is organized as follows. In the next section we introduce the
needed notations and conditions. In the section 3 we describe and recall
auxiliary facts about exponential bounds for tail of maximum distribution
of random fields. \par
 In the fourth section we will formulate and prove the main result of this
paper. Further we consider as a particular case of smooth density function.\par

 In the six  section we consider as an application the case of
{\it sample,} i.e. the case when the {\it observations}
$ \xi = \{ \xi_i, \ i = 1,2, \ldots,n \} $ are independent and identically
distributed (i., i.d). \par
 In the last section 7 we consider some examples in order to illustrate the precision of the  obtained results. \par
  Agreement: by the symbols $ C, C_j, C(i) $ we will denote some finite positive
non \ - \ essential constants.\par

\vspace{4mm}

{\bf 2. Notations and conditions. Key Inequality.}\\

\vspace{3mm}

{\sc It is presumed that all introduced function there exist in some domains which is described below. } \par

$$
U(v) \stackrel{def}{=} \{ \theta: \ \theta \in \Theta, \ r(\theta, \theta_0)
\ge v \}, \ v > 0; \eqno(2.0)
$$
then

$$
W(v) = {\bf P}(\hat{\theta} \in U(v) ). \eqno(2.1)
$$

 Let $ A(k), k = 1,2, \ldots $ be some numerical strictly increasing sequence,
$ A(1) = 1. $ For instance, $ A(k) = k $ or $ A(k) = k^{\Delta}, \
\Delta = \const > 0 $ or possible $ A(k) = C \ k^{\Delta}, \ k \ge k_0. $
 We introduce also the following measurable sets:

$$
U_k = U_k(v) = \{ \theta: \ \theta \in \Theta, \ r(\theta, \theta_0)
 \in [ \ A(k) \ v, \ A(k+1) \ v \ ] \ \}.
$$

We observe:

$$
W(v) \le \sum_{k=1}^{\infty} W_k(v), \eqno(2.2)
$$
where

$$
W_k(v) = {\bf P} \left(\hat{\theta} \in U_k(v) \right) =
{\bf P} \left( r(\hat{\theta}, \theta) \in
\left[ \ A(k) \ v, \ A(k+1) \ v \ \right] \right).
$$

 Introduce also the Kullback \ - \ Leibler "distance", or relative entropy,
or "quasi \ - \ distance"
between the parameters $ \theta $ and $ \theta_0 $ as usually

$$
h(\theta) = h(\theta_0, \theta) = {\bf E} L(\xi, \theta) =
\int_X f(x,\theta_0) \ \log[f(x,\theta)/f(x,\theta_0) ] \ d \mu. \eqno(2.3)
$$

It is known that $ h(\theta) \ge 0 $ and
$ h(\theta) = 0 \ \Leftrightarrow \theta = \theta_0. $ \par

 We denote also

 $$
 Y(v) = \inf_{\theta \in U(v) } h(\theta) \eqno(2.4)
 $$
{\it and suppose $ Y(v) \in ( 0, \infty) $ for all sufficiently great values}
 $ v. $ \par

 Further, define the following functions (some modifications of
Hellinger's integral)

$$
\phi(\lambda) = \sup_{\theta \in U_1(v) } \left[ \ \log {\bf E}
\exp \left[ \lambda L^{(0)}(\xi, \theta) \right] \ \right] =
$$

$$
\sup_{\theta \in U_1(v) } \left[ \log {\bf E} \exp( \lambda L(\xi, \theta) )
 \ \cdot \ \exp(- \lambda h(\theta)) \right] =
$$

$$
\sup_{\theta \in U_1(v) } \left[ \int_X f^{\lambda}(x,\theta) \ f^{1 \ - \ \lambda}(x,\theta_0) \ d \mu \ \cdot \ \exp(- \lambda h(\theta)) \right]; \eqno(2.5)
$$
and

$$
\chi(\lambda) = \sup_{\theta \in U(v) } \left[ \ \log {\bf E}
\exp \left[ \lambda L^{(0)}(\xi, \theta) \right] \ \right] =
$$

$$
\sup_{\theta \in U(v) } \left[ \log {\bf E} \exp( \lambda L(\xi, \theta) )
 \ \cdot \ \exp(- \lambda h(\theta)) \right] =
$$

$$
\sup_{\theta \in U(v) } \left[ \int_X f^{\lambda}(x,\theta) \ f^{1 \ - \ \lambda}(x,\theta_0) \ d \mu \ \cdot \ \exp(- \lambda h(\theta)) \right]. \eqno(2.5a)
$$

{\it We suppose $ \phi(\lambda) < \infty $ or correspondingly
$ \chi(\lambda) < \infty $ for all values $ \lambda $ in
some interval of a view} $ (0, \lambda_0), \ \lambda_0 \in (0, \infty]: $

$$
\exists \lambda_0 \in (0,\infty], \ \forall \lambda \in (0, \lambda_0] \
\Rightarrow \phi(\lambda) < \infty, \eqno(2.6)
$$

$$
\exists \lambda_0 \in (0,\infty], \ \forall \lambda \in (0, \lambda_0] \
\Rightarrow \chi(\lambda) < \infty, \eqno(2.6a)
$$

\vspace{4mm}

{\bf Key inequality.} \\

\vspace{3mm}

Note that

$$
W(v) = {\bf P} \left( \sup_{\theta \in U(v)} L(\xi,\theta) >
\sup_{\theta \notin U(v)} L(\xi, \theta) \right).
$$

As long as $ L(\xi, \theta_0) = 0 $ and $ \theta_0 \notin U $ (and
$ \theta_0 \notin U_1) $ we conclude

$$
W(v) \le {\bf P} \left( \sup_{\theta \in U(v)} L(\xi,\theta) > 0 \right) =
$$

$$
{\bf P} \left( \sup_{\theta \in U(v)} ( L(\xi,\theta)\ - \ {\bf E} L(\xi,\theta)) > \inf_{\theta \in U(v)} h(\theta) \right) \le
$$

$$
 {\bf P} \left( \sup_{\theta \in U(v)} L^0(\xi,\theta) > Y(v) \right) =
{\bf P} \left(\sup_{\theta \in U(v)} [ \log(f(\xi,\theta)/f(\xi, \theta_0))
\ - \ h(\theta) \ ] > Y(v) \right). \eqno(2.7)
$$

 Therefore, we can use the known {\it exponentially exact}
 estimations of maximum random field distributions,
 see, for example, \cite{Bagdasarova1}, \cite{Kurbanmuradov1},
\cite{Ostrovsky1}, \cite{Ostrovsky2},\cite{Ostrovsky3},
 \cite{Ostrovsky4}, \cite{Talagrand1}, \cite{Talagrand2}, \cite{Talagrand3},
 \cite{Talagrand4} etc.\par

\vspace{4mm}

{\bf 3. Auxiliary facts.}\\

\vspace{3mm}

 Let $ (\Omega, \cal{M}, {\bf P} ) $ be again the probability space,
$ \Omega = \{\omega\}, \
T = \{t \} $ be arbitrary set, $ \xi(t), \ t \in T $ be centered:
$ {\bf E} \xi(t) = 0 $ separable random field (or process). For arbitrary
subset $ S \subset T $ we denote for the values $ u \ge u_0, \ u_0 = \const
\in (0,\infty) $

$$
Q(S,u) \stackrel{def}{=} {\bf P}( \sup_{t \in S} \xi(t) > u); \
  Q(u):= Q(T,u). \eqno(3.1)
$$

{\bf Our (local) goal in this section is description
an exponentially exact as $ u \to \infty $
 estimation for the probability $ Q(S,u), \ Q(u) $ in the so \ - \ called
natural terms. } \par

  Definitions and some important results about $ {\bf E} \sup_{t \in T}
\xi(t) $ in general, i.e. non \ - \ Gaussian case, i.e. when the random field
$ \xi(t) $ may be non \ - \ Gaussian, in the terms of majoring measures see, for example, in \cite{Bagdasarova1},
\cite{Talagrand1},\cite{Talagrand2}, \cite{Talagrand3}, \cite{Talagrand4}.
In the so-called "entropy" terms this problem was considered in
\cite{Fernique1},
 \cite{Ostrovsky1}, \cite{Ostrovsky2}, \cite{Ostrovsky3}, \cite{Ostrovsky4}
 etc. \par
 In order to formulate our result, we need to introduce some addition
notations and conditions. Let $ \phi = \phi(\lambda), \lambda \in
[0, \lambda_0), \ \lambda_0 = \const \in (0, \infty] $ be some
strictly convex taking non \ - \ negative values function, such that
$ \phi(0) = 0 $ and

$$
\lambda \in [0, 0.5 \ \lambda_0) \ \Rightarrow C_1 \lambda^2 \le
\phi(\lambda) \le C_2 \lambda^2; \eqno(3.2)
$$

$$
 \lim_{\lambda \to \lambda_0} \phi(\lambda)/\lambda = \infty. \eqno(3.3)
$$

 Note that under the assumptions (2.5) or (2.5a)
$ \phi(\cdot) \in \Phi, \ \chi(\cdot) \in \Phi. $ \par

 We denote the set of all these function as $ \Phi; \ \Phi =
\{ \phi(\cdot) \}. $ \par

{\it Further we will choose $ T = U(v) $ or $ T = U_k(v) $ and $ \phi(\lambda) $
or correspondingly $ \chi(\lambda) $ as it is defined as in (2.5) and
(2.5a) }.\par

 We say that the {\it centered} random variable (r.v) $ \xi = \xi(\omega) $
belongs to the space $ B(\phi), $ if there exists some non \ - \ negative constant $ \tau \ge 0 $ such that

$$
\forall \lambda \in [0, \lambda_0) \ \Rightarrow
{\bf E} \exp(\lambda \xi) \le \exp[ \phi(\lambda \ \tau) ], \eqno(3.4)
$$
(the concretization of right \ \ hand side Kramer's condition).\par
 The minimal value $ \tau $ satisfying (3.4) is called a $ B(\phi) \ $ norm
of the variable $ \xi, $ write $ ||\xi|| $ or more detail,
$ ||\xi||B(\phi): $

 $$
||\xi|| = ||\xi||B(\phi): = \inf \{ \tau, \ \tau > 0, \ \forall \lambda \ \Rightarrow {\bf E}\exp(\lambda \xi) \le \exp(\phi(\lambda \ \tau)) \}.
\eqno(3.5)
 $$

 This spaces are very convenient for the investigation of the r.v. having a
exponential decreasing right side
 tail of distribution, for instance, for investigation of the limit theorem,
the exponential bounds of distribution for sums of random variables,
non-asymptotical properties, problem of continuous of random fields,
study of Central Limit Theorem in the Banach space etc.; see
\cite{Ostrovsky2}. \par

 The space $ B(\phi) $ with respect to the norm $ || \cdot ||B(\phi) $ and
ordinary operations is a quasi \ - \ Banach space. This means by definition
that:\\

\vspace{2mm}

{\bf 1.} $ B(\phi) $ is complete relative the quasi \ - \ distance
$ ||\xi \ - \ \eta||; $ \\

$$
{\bf 2.} ||\xi|| \ge 0; \ ||\xi|| = 0 \Leftrightarrow \xi = 0 \ -
$$
the non \ -\ negativeness;\\

$$
{\bf 3.} ||\xi + \eta|| \le ||\xi|| + ||\eta|| \ - 	
$$
the triangle inequality;

$$
{\bf 4.} \alpha = \const > 0 \ \Rightarrow ||\alpha \xi|| = \alpha ||\xi|| \ -
$$
the positive homogeneous. \par
  The $ B(\phi) $ is isomorphic to the subspace
consisted on all the centered variables of quasi \ - \
Orlitz's space $ (\Omega,F,{\bf P}), N(\cdot) $ with $ N \ - \ $ right function

$$
N(u) = \exp(\phi^*(u)) - 1, \ \phi^*(u) = \sup_{\lambda} (\lambda u -
\phi(\lambda)).
$$
 The transform $ \phi \to \phi^* $ is called Young \ - \ Fenchel or Legendre
 transform. The proof of considered assertion used the properties of saddle-point method and theorem of Fenchel \ - \ Moraux:
$$
\phi^{**} = \phi.
$$

 Many facts about the $ B(\phi) $ spaces are proved in \cite{Ostrovsky2},
\cite{Ostrovsky3}, p. 19 \ - \ 40:

$$
 \xi \in B(\phi) \Leftrightarrow {\bf E } \xi = 0, \ {\bf and} \ \exists C = \const > 0,
$$

$$
Z(\xi,x) \le \exp \left(-\phi^*(Cx) \right), x \ge 0, \eqno(3.6)
$$
where $ Z(\xi,x)$ denotes in this article the {\it right hand tail} of
distribution of the r.v. $ \xi: $

$$
Z(\xi,x) \stackrel{def}{=} {\bf P}(\xi > x), \ x \ge 0,
$$
and this estimation is in general case asymptotically exact. \par
More exactly, if $ \lambda_0 = \infty, $ then the following implication holds:

$$
\lim_{\lambda \to \infty} \phi^{-1}(\log {\bf E} \exp(\lambda \xi))/\lambda =
K \in (0, \infty)
$$
if and only if

$$
\lim_{x \to \infty} (\phi^*)^{-1}( |\log Z(\xi,x)| )/x = 1/K.
$$
 Here and further $ f^{-1}(\cdot) $ denotes the inverse function to the
function $ f $ on the left \ - \ side half \ - \ line $ (C, \infty). $ \par

 Let $ \phi \in \Phi. $ We denote

$$
\phi_n(\lambda) = n \phi(\lambda/\sqrt{n}), \
\overline{\phi}(\lambda) = \sup_{n=1, 2, \ldots} [n \ \phi(\lambda/\sqrt{n})]
\eqno(3.7)
$$
and analogously

$$
\chi_n(\lambda) = n \chi(\lambda/\sqrt{n}), \
\overline{\chi}(\lambda) = \sup_{n=1, 2, \ldots} [n \ \chi(\lambda/\sqrt{n})].
\eqno(3.7a)
$$

 This function obeys the following sense. If $ \xi(i), \ i = 1,2,\ldots $ be
a sequence of centered, i., i.d. r.v., $ \xi = \xi(1), $ belonging to the
space $ B(\phi) $ and having the unit norm in this space: $ ||\xi||B(\phi) 
 = 1, $ then we have for the normed sum

$$
\eta(n) = n^{-1/2} \ \sum_{i=1}^n \xi(i):
$$

$$
{\bf E} [\exp(\lambda \ \eta(n) ) ] \le \exp[\phi_n(\lambda) ],
$$

$$
\sup_{n=1,2,\ldots} {\bf E} \exp(\lambda \ \eta(n)) \le
\exp \left[\overline{\phi}(\lambda) \right], \eqno(3.8)
$$
and following

$$
 Z(\eta(n),x) \le \exp \left[ \ - \ (\phi_n)^*(x) \right],
$$
the non \ - \ uniform estimation;

$$
\sup_n Z(\eta(n),x) \le \exp \left[ \ - \ (\overline{\phi})^*(x) \right],
\eqno(3.9)
$$
the uniform estimation (Chernoff's estimations, see \cite{Chernoff1}). \par
 Using the property (3.2), we can show in addition to the classical theory of
 the great deviations that in the "mild" zone

$$
x = x(n) \in (0, C \ \sqrt{n}) \ \Rightarrow
$$

$$
Z(\eta(n), x) \le C_2 \ \exp \left( - C_3 \ x^2 \right)
$$
(the non \ - \ uniform estimation). \par

As an example: if in addition

$$
Z(\xi(i),x) \le \exp \left( - x^q \right), \ q = \const \ge 1, \
x \ge 0, \eqno(3.10)
$$
then for some "constant" $ C = C(q) \in (0,\infty) $

$$
\sup_{n = 1,2,\ldots} Z(\eta(n), x) \le \exp \left[- C x^{\min(q,2)} \right],
\eqno(3.11)
$$	
and the last estimation is unimprovable at $ x >>1. $\par

 Now we prove a more general assertion.\par

\vspace{2mm}

{\bf Lemma 3.1} Let $ \{ \eta(i) \}, \ i = 1,2,\ldots, n $ be a sequence of
i., i.d., centered r.v. such that for some $ q = \const > 0 $ and for all positive values $ x $

$$
Z(|\eta(i)|, x) \le \exp \left( - x^q \ R(x) \right), \eqno(3.12)
$$
where $ R(x) $ is continuous positive {\it slowly varying} as $ x \to \infty: $

$$
\forall t > 0 \ \Rightarrow \lim_{x \to \infty} R(tx)/R(x)  = 1;
$$
is bounded from below in the  positive semi \ - \  axis

$$
 \inf_{x \ge 0} R(x) > 0
$$
function. For instance: $ R(x) = $

$$
 [\log(x + 3)]^r,   \ R(x) = [\log(x + 3)]^r  \cdot
[\log(\log(x + 16))]^s, \ r,s = \const, r \ge 0.
$$
 Denote

 $$
 \zeta(n) = n^{-1/2} \sum_{i=1}^n \eta(i).
 $$
We assert:

$$
\sup_n Z(|\zeta(n)|, x) \le
\min \left[ \exp \left( - C_1(q,R) \ x^q \ R(x) \right), \
\exp \left( - C_2(q,R) \ x^2 \right) \right]. \eqno(3.13)
$$

Notice that the lower bound, i.e. the {\it inverse} inequality

$$
\sup_n Z(|\zeta(n)|, x) \ge
\min \left[ \exp \left( - C_3(q,R) \ x^q \ R(x) \right), \ \exp \left( - C_4(q,R) \ x^2 \right) \right] \eqno(3.14)
$$
in the case when

$$
Z(|\eta(i)|,x) \le \exp \left[ - C_5(q,R) \ x^q \ R(x) \right]
$$
and
$$
Z(|\eta(i)|,x) \ge \exp \left[ - C_6(q,R) \ x^q \ R(x) \right], \ 0 < C_5 \ge C_6 < \infty
$$
is trivial. Namely,

$$
\sup_n Z(|\zeta(n)|, x) \ge Z(|\eta(1)|,x) \ge \exp \left[ - C_6(q,R) \ x^q \ R(x) \right],
$$
and on the other hand

$$
\sup_n Z(|\zeta(n)|, x) \ge \lim_{n \to \infty} {\bf P}(|\zeta(n)| > x ) =
$$

$$
2 \int_x^{\infty} (2 \pi)^{-1/2} \ \sigma^{-1} \ \exp
\left[-y^2/(2 \sigma^2) \right] \ dy \ge
\exp \left[ - C_7(q,R) \ x^2 \right], \ x \ge C_8 > 0;
$$
here we used the CLT and denoted
$$
\sigma^2 = \sigma^2(q) = {\bf Var} \ (\eta(1)) \in (0, \infty).
$$

\vspace{2mm}
{\bf Proof} (briefly) of the Lemma 3.1. \par

 The case $ q \ge 1 $ is considered in \cite{Ostrovsky2}, chapter 1,
section (1.6); therefore we must consider only the case $ q \in (0,1). $ \par

  Further, without loss of generality we can consider the case when the r.v.
 $ \xi $ and $ \eta $ are independent and symmetrical distributed with densities
 correspondingly
 $$
 f_{\xi}(x) = f(x) = C_9 \ \exp \left( - |x|^q \ R(|x|) \right),
 $$

 $$
 g_{\eta}(x) = g(x) = C_{10} \ \exp \left( - K \ |x|^q \ R(|x|) \right),
 $$
 and $ \tau = \xi + \eta. $ Here
 $ K = \const \in (1, \infty) \ $ (the case when $ K = 1 $ may be considered analogously). \par
 Let us assume that $ x \to \infty, \ x \ge C. $ We have denoting by $ h(x) =
h_{\tau}(x) $ the density of distribution of the r.v. $ \tau: $

$$
h(x) \sim C_{11} \int_0^x \exp \left[ - K(x \ - \ y)^q R(x \ - \ y) \ - \
y^q R(y) \right] \ dy =
$$

$$
C_{11} x \int_0^1 \exp \left[ - x^q \left[ K (1 \ - \ t)^q R(x \ (1 \ - \ t)) +
t^q \ R(t \ x)  \right]    \right] \  dt \sim
$$

$$
C_{11} x \int_0^1 \exp \left[ - x^q \ R(x) \ S(t) \right] \ dt,
$$
where
$$
S(t) = K(1 \ - \ t)^q  + t^q.
$$
  The function $ t \to S(t), \ t \in [0,1] $ achieves the minimal value $ K $
at the (critical) point $ t = 0 $ and as $ t \to 0+ $

$$
S(t) = K + t^q + 0(t).
$$

 Note that in the case $ K = 1 $ there are two critical points:
 $ t = 0 $ and  $ t = 1.$ \par

Further, we use the classical saddle \ - \ point method: at $ x \to \infty, x > 1  $ we have:

$$
h(x) \sim C_{11} x \int_0^{\infty} \exp \left[ - x^q \ R(x) \ (K + t^q) \right] \ dt =
$$

$$
C_{12}(q, R) \ x \ \exp \left[ - K \ x^q \ R(x)  \right] \ \left(x^q \ R(x)
\right)^{-1/q} \le
$$

$$
C_{13}(q,R) \ \exp \left[ - K \ x^q \ R(x)  \right].
$$

 This completes the proof of the lemma 3.1. \par

\vspace{4mm}

The function $ \phi(\cdot) $ may be introduced  "constructive", i.e. {\it only
by means of the values of the considered random field} 
$ \{ \xi(t), \ t \in T \} $  by the formula

$$
\phi(\lambda) = \phi_0(\lambda) \stackrel{def}{=} \log \sup_{t \in T}
 {\bf E} \exp(\lambda \xi(t)), \eqno(3.15)
$$
 if obviously the family of the centered r.v. $ \{ \xi(t), \ t \in T \} $ satisfies the {\it uniform } Kramer's condition:
$$
\exists C \in (0, \infty), \ \sup_{t \in T} Z(\xi(t), \ x)
\le \exp(- C \ x), \ x \ge 0. \eqno(3.16)
$$
 In this case, i.e. in the case the choice the function $ \phi(\cdot) $ by the
formula (3.15), we will call the function $ \phi(\lambda) = \phi_0(\lambda) $
a {\it natural } function. \par
 Note that if for some $ C = \const \in (0,\infty) $
$$
Q(T,u) \le \exp \left(-\phi^*(C u) \right),
$$
 then the condition (2.6) is satisfied (the "necessity" of the condition (2.6)). \par
 M.Talagrand \cite{Talagrand1}, \cite{Talagrand2},
 \cite{Talagrand3}, \cite{Talagrand4},
W.Bednorz \cite{Bednorz1}, X. Fernique \cite{Fernique1}
 etc. write
instead our function $ \exp \left(-\phi^*(x) \right) $  the function 
$1/\Psi(x), $ where $ \Psi(\cdot) $ is some Young's
function and used as a rule a function $ \Psi(x) = \exp(x^2/2) $
(the so \ - \ called "subgaussian case").\par

 Without loss of generality we can and will suppose

$$
\sup_{t \in T} [ \ ||\xi(t) \ ||B(\phi)] = 1,
$$
(this condition is satisfied automatically in the case of natural choosing
of the function $ \phi: \ \phi(\lambda) = \phi_0(\lambda) \ ) $
and that the metric space $ (T,d) $ relatively the so-called
{\it natural} distances (more exactly, semi \ - \ distances)

$$
d_{\phi}(t,s) = d(t,s) \stackrel{def}{=} ||\xi(t) \ - \ \xi(s)|| B(\phi) ,\eqno(3.17)
$$
and analogously  (see the definition of a function $ \chi(\cdot) $ further)\par

$$
d_{\chi}(t,s) \stackrel{def}{=} ||\xi(t) \ - \ \xi(s)|| B(\chi) \eqno(3.17a)
$$
is complete. \par
 Recall that the semi \ - \ distance $ \rho = \rho(t,s), \ s,t \in T, $
for instance, $ d = d_{\phi}(t,s), \ s,t \in T $ is, by definition, a
non \ - \ negative symmetrical numerical function, $ \rho(t,t) = 0, \ t \in T, $
satisfying the triangle inequality, but the equality $ \rho(t,s) = 0 $
does not means (in general case) that $ s = t. $ \par
 For example, if the random field $ \xi(t) $ is centered and normed:

$$
\sup_{t \in T} {\bf Var} \ [\xi(t)] = 1
$$
Gaussian field with a covariation function
 $ D(t,s) = {\bf E} \xi(t) \ \xi(s), $ then
$ \phi_0(\lambda) = 0.5 \ \lambda^2, \ \lambda \in R, $ and $ d(t,s) = $

$$
||\xi(t) - \xi(s)||B(\phi_0) = \sqrt{ \bf{Var} [ \xi(t) - \xi(s) ]} =
\sqrt{ D(t,t) \ - \ 2 D(t,s) + D(s,s) }.
$$

 Let $ (T, \rho) $ be a compact metrical space.  We us introduce as usually
 for any subset $ V, \ V \subset T $ the so-called
{\it entropy } $ H(V, \rho, \epsilon) = H(\rho, \epsilon) $ as a logarithm
of a minimal quantity $ N(V,\rho, \epsilon) = N(V,\epsilon) = N(\rho,
\epsilon)  $ of a balls in the distance $ \rho(\cdot, \cdot) $
$ S(V, t, \epsilon), \ t \in V: $
$$
S(V, t, \epsilon) \stackrel{def}{=} \{s, s \in V, \ \rho(s,t) \le \epsilon \},
$$
which cover the set $ V: $
$$
N = \min \{M: \exists \{t_i \}, i = 1,2,…, M, \ t_i \in V, \ V
\subset \cup_{i=1}^M S(V, t_i, \epsilon ) \},
$$
and we denote also

$$
H(V,\rho,\epsilon) = \log N; \ S(t_0,\epsilon) \stackrel{def}{=}
 S(V, t_0, \epsilon),
$$

$$
 H(\rho, \epsilon)  \stackrel{def}{=} H(T,\rho,\epsilon). \eqno(3.18)
$$
 It follows from Hausdorf's theorem conversely that
$ \forall \epsilon > 0 \ \Rightarrow H(V,\rho,\epsilon)< \infty $ iff the
metric space $ (V, \rho) $ is precompact set, i.e. is the bounded set with
compact closure.\par

 We quote now some results from \cite{Ostrovsky2}, \cite{Ostrovsky3},
\cite{Ostrovsky4} about the non \ - \ asymptotic 
exponential estimations for $ Q(u) = Q(T,u) $ as
$ u >> 1. $ Define for any value $ \delta \in (0, 1) $ and arbitrary subset
$ V $ of the space $ \Theta: \ V \subset \Theta $ and some semi \ - \ distance
$ \rho(\cdot, \cdot) $  on the set $ T $  the following function:

$$
  G(V,\rho,\delta) = G(\rho, \delta) =
 \sum_{m=1}^{\infty} \delta^{m \ - \ 1} \ \cdot \ H(V, \rho, \delta^m) \ \cdot
 \ (1 \ - \ \delta). \eqno(3.19)
$$
 We define formally $ G(\delta) = + \infty $ for the values $ \delta
> \delta_0. $ \par
 In the case when $ V = U(v) $ and $ \rho(t,s) = d_{\phi}(t,s), $
i.e.  when $ \rho $ is the natural semi \ - \ distance,
we will write  for brevity $ G(\delta) = G(U(v), d_{\phi}, \delta). $ \par

 If

$$
\exists \delta_0 \in (0,1), \ \forall \delta \in (0, \delta_0) \
\Rightarrow G(\delta) < \infty, \eqno(3.20)
$$
then

$$
Q(T,u) \le V(T, \delta,u), \
V(T, \delta,u) \stackrel{def}{=} \exp \left( G(\delta) \ - \
\phi^*(u(1 \ - \ \delta)) \right), \eqno(3.21)
$$
or equally

$$
 Q(T,u) \le \inf_{\delta \in (0,1)} V(T, \delta,u). \eqno(3.22)
$$

 If for example

$$
\forall \delta \in (0, 1/e] \ \Rightarrow
G(T, d_{\phi}, \delta) \le H_0 + \kappa |\log \delta|, \ H_0, \kappa = \const < \infty,
$$
then we get denoting

$$
\pi(u) = u \ \phi^{*/ }(u)	
$$
for the values $u $ for which $ \pi(u) \ge 2 \kappa: $

$$
Q(T,u) \le \exp(H_0) \ C^{\kappa} \ \kappa^{-\kappa} \ (\pi(u))^{\kappa} \
\exp \left(-\phi^*(u) \right), \eqno(3.23)
$$
and the last estimation (3.23) is exact in the main ("exponential") term
$ \exp(\left(-\phi^*(u) \right). $ \par
 More exactly, in many practical cases the following inequality holds:

$$
\forall \epsilon \in (0,3/4) \ \exists K > 0, \ \forall u > K \ \Rightarrow \pi(u)
< \exp( \ - \ \phi^*(\epsilon \ u)); \eqno(3.24)
$$
and we conclude  hence for  $ u > K = K(\epsilon) $ by virtue of convexity of a function $ \phi^*(x): $

$$
Q(T,u) \le \ C_1(\kappa, \phi(\cdot) ) \
\exp \left(-\phi^*( (1 \ - \epsilon) \ u) \right),
$$
and conversely there exists a r.v. $ \zeta $ with unit norm in the space
$ B(\phi): \ \zeta \in B(\phi), \ ||\zeta|| = 1,$
for which

$$
u \ge K \ \Rightarrow \ Z(\zeta,u) \ge \ C_2(\phi) \
\exp \left(-\phi^*( (1 \ + \epsilon) \ u) \right).
$$

 The value $ \kappa $ is called the {\it metric dimension} of the set
$ T $ relative the distance $ d = d_{\phi}(\cdot, \cdot). $ \par

 Note that if
$$
T = \cup_{m=1}^{\infty} T(m)
$$
is some measurable partition $ R = \{ T(m) \} $ of the parametrical set $ T, $ then

$$
Q(T,u) \le \sum_{m=1}^{\infty} Q(T_m,u)
$$
and hence

$$
Q(T,u) \le \inf_{ R = \{T(m) \}} \sum_{m=1}^{\infty} Q(T_m,u).
$$
 Estimating the right side term by means of the inequality (3.19), we get:

$$
Q(T,u) \le \inf_{ R = \{T(m) \}}
\left[ \sum_{m=1}^{\infty} \ \inf_{\delta(m) \in (0,1)}
 \sum_{m=1}^{\infty} V(T(m), \delta(m),u) \right]. \eqno(3.24)
$$

 The last assertion is some simplification of the {\it Majorizing Measures,
or Generic Chaining Method} (see \cite{Fernique1},
\cite{Talagrand1} \ - \ \cite{Talagrand4}, \cite{Bednorz1}, \cite{Ledoux1}
 etc). \par

Further we will use as a rule the partition $ R $  of the set $ U(v) $ of a
view

$$
R = \cup_{k=1}^{\infty}  \{ \theta: \theta \in U(v), \ r(\theta, \theta_0)
\in [ \ A(k) \ v, \ A(k+1) \ v  \ ] \ \}. \eqno(3.25)
$$

\vspace{4mm}

{\bf 4. Main results.} \\

\vspace{3mm}

{\bf A. "Compact" parametrical set.} \\

\vspace{2mm}

The "compactness" means by definition that the function $ \theta \to
r(\theta, \theta_0), \theta \in \Theta $ is bounded. Since as a rule 
the parametric set
$ \Theta $ is a closed subset in Euclidean finite \ - \ dimensional space
and $ r(\cdot, \cdot) $ is ordinary distance, this definition coincides
with usually definition of the compact sets. \par
 Note that in this case only finite numbers of the sets $ \{ A(k) \} $ are
non \ - \ empty. We can suppose in this subsection for simplicity 
$ U_1(v) = U(v) $ and therefore $ \phi(\lambda) = \chi(\lambda). $ \par

 Let the function $ \phi = \phi(\lambda) $ be defined as in (2.5a) or equally
(in the considered case) as in (2.5). Recall that
$$
\sup_{\theta \in U(v)} || L^0(\xi, \theta)||B(\chi) = 1.
$$

 Introduce the so \ - \ called {\it natural} semi \ - \ distance on the set
$ U(v) $ as follows:

$$
 d = d_{\chi} = d(\theta_1, \theta_2) = ||L^0( \xi, \theta_1) \ - \
L^0( \xi, \theta_2)||B(\chi)  =
$$

$$
|| \log [f(\xi, \theta_1)/f(\xi, \theta_2) ]  \ - \ [h(\theta_1) \ - \ h(\theta_2)]||B(\chi). \eqno(4.1)
$$

 It follows immediately from (3.18) (or equally from (3.19)) the following result. \par

{\bf Theorem 4.1.a.}. If there exists $ \delta_0 = \const \in (0,1) $ such that
$ \forall \delta \in (0,1) \ \Rightarrow $

$$
G(U(v), d_{\chi}, \delta) := \sum_{m=1}^{\infty} \delta^{m \ - \ 1} \ H(U(v), d_{\chi}, \delta^m) \ (1 \ - \ \delta) < \infty, \eqno(4.2),
$$
then $ \forall \delta \in (0, \delta_0] $

$$
W(v) \le \exp \left[ G( U(v), d_{\chi}, \delta) \ - \ \chi^*( (1 \ - \ \delta) \ Y(v) ) \right]. \eqno(4.3)
$$

 Let us offer the more convenient for application form. Define for
$ \tilde{U} \subset \ U(v),$  arbitrary function $ \nu \in \Phi, $ and any
semi \ - \ distance  $     \rho = \rho(\theta_1, \theta_2) $ on the set
$ \tilde{U} $ the following  function (if it is finite)

$$
\Psi_{\nu}( \tilde{U}, \rho, y) \stackrel{def}{=} \inf_{\delta \in (0,1)}
\exp \left[ G(\tilde{U}, \rho, \delta) \ - \ \nu^*((1 \ - \ \delta) \ y )  \  \right].\eqno(4.4)
$$

{\bf Theorem 4.1.} Under the conditions of the theorem (4.1.a) the following
estimate is true:

$$
W(v) \le \Psi_{\chi}(U(v), d_{\chi}, Y(v)). \eqno(4.5)
$$

\vspace{3mm}

{\bf B. "Non \ - \ compact" set.}\\

\vspace{3mm}

 In this case we need to use the main idea of the so \ - \ called "generic
chaining ", or "majorizing measure" method (3.22), (see \cite{Fernique1},
\cite{Talagrand1} \ - \ \cite{Talagrand4}, \cite{Bednorz1}, \cite{Ledoux1}
 etc), which used in particular the partition  $ U(v)  = \cup_k U_k(v). $ \par

Let us denote for the partition $ R = \{ U_k(v) \}, U(v) = \cup_k U_k(v) $

$$
\tau_k = \tau_k(v) = \sup_{\theta \in U_k(v)} ||L^0(\xi, \theta)||B(\phi),
$$

$$
Y_k(v) = \inf_{\theta \in U_k(v)} h(\theta),
$$
and introduce the following distance $ d_k $ on the set $ U_k = U_k(v): $

$$
d_k \left(\theta_1^{(k)}, \theta_2^{(k)} \right) =
||L^0(\xi, \theta_1^{(k)}) \ - \ L^0(\xi, \theta_2^{(k)}) ||B(\phi),
 \theta_1^{(k)}, \  \theta_2^{(k)} \in U_k(v).
$$

{\bf Theorem 4.2}. We have for arbitrary partition $ R $

$$
W(v) \le \sum_{k=1}^{\infty} \Psi_{\phi}
\left(U_k(v), \frac{d_k}{\tau_k(v)}, \frac{Y_k(v)}{\tau_k(v)} \right).
\eqno(4.6)
$$

Notice that

$$
H(V, \rho/K, \epsilon) = H(V, \rho, K \cdot \epsilon), \ K = \const > 0.
\eqno(4.7).
$$

{\bf Proof} of the Theorem 4.2.  We use the inequality (2.2): $ W(v) \le
\sum_k W_k(v).$  Let us estimate each summand $ W_k(v): $

$$
W_k(v) = {\bf P} \left(\sup_{\theta \in U_k(v) } L(\xi, \theta) > 0 
\right) = 
$$

$$
{\bf P} \left(\sup_{\theta \in U_k(v) } [ L^0(\xi, \theta) \ - \
h(\theta)] > 0 \right) \le
$$

$$
{\bf P} \left( \sup_{\theta \in U_k(v) } L^0(\xi, \theta) > Y_k(v)
 \right) = {\bf P} \left( \sup_{\theta \in U_k(v) }
\frac{ L^0(\xi, \theta)}{\tau_k(v) }  >  \frac{ Y_k(v) }{\tau_k(v) } \right).
\eqno(4.8)
$$

 The random field

$$
\xi_k(\theta) = \frac{ L^0(\xi, \theta)}{\tau_k(v) }, \ \theta \in U_k(v)
$$
is normed in the $ B(\phi) $ sense:

$$
\sup_{\theta \in U_k(v)} ||\xi_k(\theta)||B(\phi) = 1.
$$

Further,

$$
|| \xi_k(\theta_1^{(k)}) \ - \ \xi_k(\theta_2^{(k)})||B(\phi) =
d_k \left(\theta_1^{(k)}, \theta_2^{(k)} \right).
$$

Using the inequality (3.22) for the probability $ W_k(v) $ and summing over
$ k, $ we arrive to the estimation (4.6).\par

\vspace{4mm}

{\bf 5. The regular, or "smooth" case.} \\

\vspace{3mm}

{\bf A. Non \ - \ formal introduction. Restrictions. Conditions.} \par

\vspace{2mm}

 In this section we consider the case when the set $ \Theta $ is closed (may be
unbounded) convex nonempty subset of the Euclidean space $ R^m, \ m = 1,2,\ldots, $ the density $ f(x,\theta) $ is twice differentiable function on
the variable (variables) $ \theta. $  \par

 We choose as the deviation function {\it hereafter} $ r(\theta, \theta_0) $
the ordinary Euclidean distance

$$
r(\theta_1, \theta_2) = \sqrt{ (\theta_1 \ - \ \theta_2, \theta_1 \ - \
\theta_2) }  \stackrel{def}{=} |\theta_1 \ - \theta_2|.
$$
 The function $ \phi $ is  in this section the natural, i.e.
$ \phi(\lambda) = \phi_0(\lambda). $ \par

We have formally  as
$ \theta \to \theta_0, $ denoting $ \nabla \ f = \grad_{\theta} \ f =
\partial f/\partial \theta: $

$$
h(\theta) \sim  \int_X f(x,\theta_0) \ \times
$$

$$
\log  \left( \frac{f(x,\theta_0) + \nabla f(x,\theta_0) (\theta \ - \
\theta_0) + 0.5 \nabla^2 \ f(x,\theta_0) (\theta - \theta_0, \theta \ - \ \theta_0)} {f(x, \theta_0) } \right) \ \mu(dx) \asymp
$$

$$
C \ r(\theta, \theta_0)^2  = C \ |\theta \ - \ \theta_0|^2, \ C = C(f(\cdot, \cdot), \theta_0).
$$

It is reasonable to assume that

$$
h(\theta) \asymp C \ |\theta \ - \ \theta_0|^2. \eqno(5.1)
$$

\vspace{3mm}

{\bf B.  Main result of this section.}

\vspace{3mm}

{\bf Theorem 5.1.} \par

\vspace{3mm}

{\it We impose on the our statistical structure the following
conditions. } \par

{\bf A.} \ Let the function
  $ \phi(\lambda) = \phi_0(\lambda) $ satisfied the
condition (2.6) on the set $ T = U(1).$ \par

{\bf B.} \  Assume that the condition (5.1) is satisfied.\par

{\bf C.} \ Suppose there exists a constant $ C > 1 $ such that
for each constant $ K > 1 $ the following inequality holds:

$$
\sup_{\theta: v \le |\theta \ - \ \theta_0| \le K \ v }
\frac{||L^0(\xi, \theta) ||} { |\theta \ - \ \theta_0 |} \le C \cdot K; \eqno(5.2)
$$

{\it  Then there exists a constant $ C = C( f(\cdot,\cdot),m, \theta_0)  \in
(0, \infty ) $ such that for all the values } $ v \ge 1 $

$$
W(v) \le \exp \left( - \phi^*(C \cdot v) \right). \eqno(5.3)
$$

\vspace{2mm}

{\bf Proof.} \par
{\bf 1.}  We intend to use the result of the theorem 4.2. First of all
we choose the partition $ R $ of a view: $ R = \cup_k [ \ A(k) \ v,
\ A(k + 1) \ v \ ], $ where $ A(k) = k, \ k = 1,2,\ldots. $ \par

\vspace{2mm}

{\bf 2.} From the conditions {\bf B}, or equally the condition (5.1)
and  the condition {\bf C }  follows that:

$$
\tau_k(v) \le C_2 \ (k + 1) \ v \eqno(5.4)
$$
and

$$
Y_k \stackrel{def}{=} Y_k(v) \ge C_3 A_k^2 v^2. \eqno(5.5)
$$

\vspace{2mm}

{\bf 3.} Since the function $ \phi(\cdot) = \phi_0(\cdot) $ satisfies the condition {\bf A }, we can estimate the "natural" distance $ d_k $ as
follows:

$$
d_k \left(\theta_1^{(k)}, \theta_2^{(k)} \right)/\tau_k =
||L^0(\xi, \theta_1^{(k)}) \ - \ L^0(\xi, \theta_2^{(k)}) ||B(\phi)
/ \tau_k \le
$$

$$
C_4 \ | \theta_1^{(k)} \ - \  \theta_2^{(k)} |. \eqno(5.6)
$$

 Since the "layer" $ U_k $ is bounded in the Euclidean metric,  we conclude
from  (5.6) that

$$
H(U_k(v),d_k/\tau_k, \delta) \le C_6(L,m) + m \ |\log \delta|. \eqno(5.7)
$$
 On the other words, in the considered  "regular" case $ \kappa = m. $ \par
Therefore, all the conditions of theorem 4.2 are satisfied, and we obtain
from the inequality (4.6): $ W(v) \le W_0(v), $ where

$$
W_0(v) \stackrel{def}{=} \sum_{k=1}^{\infty}
\exp \left( - \phi^*(k^2 \ v^2/(C_7  \ (k+1) \ v  ) ) \right) \le
$$

$$
\sum_{k=1}^{\infty} \exp \left( - \phi^*(C_8 \ k \ v)  \right) \ \le \
\exp \left( - \phi^*( C_9 \ v ) \right) \eqno(5.8)
$$
as long as $ v \ge 1.$ \par
 This completes the proof of theorem 5.1. \par

\vspace{3mm}

{\bf Corollary 5.1.} \
The conclusion of the theorem (5.1), i.e. the inequality (5.3) may be
rewritten  as follows. For all the values $ v \ge 0 $

$$
W(v) \le \min \left[1, W_0(v)  \right]. \eqno(5.9)
$$

 Note that $ W_0(0) = + \infty. $ \par

\vspace{3mm}

{\bf Corollary 5.2.} \

We obtain using the asymptotical behavior of the function $ \phi = \phi(\lambda),
\lambda \to 0+ $ in the {\it bounded interval } of the variable
$ \ v: \  v \in [1, C_1], \ C_1 = \const > 1  $

$$
W(v) \le \exp \left( - C \ v^2 \right). \eqno(5.10)
$$

Notice that under some additional conditions, see \cite{Ibragimov1}, chapter 3, section 3, at $ v \le 1 $ the following inequality holds:

$$
W(v) \le \exp \left( - C_2 \ v^2 \right).
$$

Therefore, we get under these conditions at $ v \le C_3 = \const > 0 $

$$
W(v) \le \exp \left( - C_4 \ v^2 \right). \eqno(5.11)
$$

\vspace{3mm}

{\bf Remark 5.1} We conclude in the "smooth" case,
taking the union of inequalities 5.8 and 5.11 and taking into account the behavior of the function $ \phi^*(\lambda) $ as $ \lambda \to 0+: $ as in the
case of the of the function $ \phi(\lambda): $

$$
\phi^*(\lambda) \sim \lambda^2, \ \lambda \in [0, C],
$$
that

$$
W(v) \le \exp \left[ - \phi^*( C \ v) \right]. \eqno(5.12)
$$

 We obtained the main result of this report.\par

\vspace{3mm}

{\bf 6. The case of sample.} \\

\vspace{3mm}

 In this section we consider the case when $ \xi = \vec{\xi} = \{\xi(i)\},
i = 1,2, \ldots,n $ are i., i.d. r.v. with the ("one \ - \ dimensional")
density $ f(x, \theta), $ {\it satisfying all the condition of the sections
1 and 5}, ( the "smooth case".) \par
 We keep also all notations for the function $ f(\cdot, \cdot), $ for instance
the notions  $ h(\theta), Y(v), \phi(\cdot), \overline{\phi}, \ R  $ etc.\par

 We will investigate in this section the {\it non \ - \ uniform} probability
{\it under natural norming } $  \sqrt{n}: $
$$
W_n(v) = {\bf P} (\sqrt{n} \ r(\hat{\theta_n}, \theta_0 ) > v), \eqno(6.1)
$$
where $ \hat{\theta}_n = \hat{\theta} $ is the MLE estimation of the unknown
parameter $ \theta_0 $ on the basis the sample $ \xi = \vec{\xi}: $

$$
\hat{\theta} = \hat{\theta_n} = \argmax_{\theta \in \Theta}
\prod_{i=1}^n f(\xi(i), \theta)
$$
or equally

$$
\hat{\theta} = \argmax_{\theta \in \Theta} L(\xi, \theta) =
\argmin_{\theta \in \Theta}(- L(\xi, \theta)) \eqno(6.2)
$$
where the contrast function $ L(\cdot, \cdot) $ may be written here as

$$
 L = L^{(n)} = L^{(n)}(\xi, \theta) = \sum_{i=1}^n
\log \left[ f(\xi(i), \theta)/f(\xi(i), \theta_0) \right] \eqno(6.3)
$$
and correspondingly 

$$
 L_0 = L^{(n)}_0 = L^{(n)}(\xi, \theta) = \sum_{i=1}^n
\log \left[ f(\xi(i), \theta)/f(\xi(i), \theta_0) \ - \ h(\theta) \right] \eqno(6.3.a)
$$

and find also the upper estimation for the {\it uniform} probability

$$
\overline{W}(v) = \sup_n W_n(v). \eqno(6.4)
$$

\vspace{3mm}

{\bf Theorem 6.1} Under the formulated conditions the following estimations
are true:

$$
W_n(v) \le \exp \left( - \phi_n^*( C_1 \ v)  \right), \eqno(6.5)
$$

$$
\overline{W}(v) \le \exp \left( -  \overline{ \phi}^*(C_1 \ v \right). \eqno(6.6)
$$

{\bf Proof.} Let us denote for brevity

$$
\eta(i, \theta) =
\eta(i) = \log \left[f(\xi(i), \theta)/f(\xi(i), \theta_0) \right],
$$

$$
\eta^o(i) = \eta^o(i,\theta)= \eta(i) \ - \ {\bf E} \eta(i) = \eta(i) \ - \ h(\theta).
$$

 We have using the key inequality for the sample of a volume $ n: $

$$
W_n(v) \le {\bf P} \left( \sup_{\theta \in U(v/\sqrt{n})}  \frac{1}{\sqrt{n}}
\sum_{i=1}^n \eta(i,\theta) > 0 \right) \le
$$

$$
{\bf P} \left( \sup_{\theta \in U(v/\sqrt{n}) }  \frac{1}{\sqrt{n} }
\sum_{i=1}^n  \left[\eta^o(i,\theta) \ - \ h(\theta) \right] > 0 \right) \le
$$

$$
{\bf P} \left( \sup_{\theta \in U(v/\sqrt{n}) }  \frac{1}{\sqrt{n}}
\left[ \sum_{i=1}^n  \eta^o(i,\theta) \right] >
 \sqrt{n} \ Y(v/\sqrt{n} ) \right)=  
$$

$$
{\bf P} \left( \sup_{\theta \in U(v/\sqrt{n}) }  \frac{1}{\sqrt{n}}
\left[ \sum_{i=1}^n  \eta^o(i,\theta) \right]/\tau(v/\sqrt{n}) >
 \sqrt{n} \ Y(v/\sqrt{n} )/ \tau(v/\sqrt{n}) \right).  \eqno(6.7)
$$

 As long as 

$$
Y(v) \ge C_1 \ v^2, \ \tau(v) \stackrel{def}{=} \tau_1(v) \le C_2 \ v, \ v > 0, 
$$

we can use for the estimation of the distribution of the r.v.

$$
\zeta_n(\theta) \stackrel{def}{=} \frac{1}{\sqrt{n}} \left[ \sum_{i=1}^n
\eta^o(i,\theta) \right] / \tau(v/\sqrt{n})
$$
and  the difference

$$
\zeta_n(\theta_1) \ - \ \zeta_n(\theta_2) =
 \frac{1}{\sqrt{n}} \left[ \sum_{i=1}^n ( \eta^o(i,\theta_1) \ - \
\eta^o(i,\theta_2) ) \right] / \tau(v/\sqrt{n}) \eqno(6.8)
$$
 the definition of the function $ \phi^*(\cdot) $ and its
properties; another approach in the many general cases, i.e. when the function
$ \phi(\cdot) $ does not exists, may be investigated by means of the Lemma 3.1. \par

Using the estimation (5.12), we affirm

$$
W_n(v) \le \exp \left( - \phi_n^* \left[ C \ v  \right] \right). \eqno(6.9)
$$

 The second assertion of the theorem 6.1 follows immediately by passing to
$ \sup_n. $ \par

\vspace{3mm}

{\bf Remark 6.1} From the assertion of the theorem 6.1 it may be obtained the
estimations from {\it integral } measures of deviation. For instance, if we choose the loss function $ l(\cdot) $ of a kind

$$
l = \sqrt{n} \ r(\hat{\theta_n}, \theta),
$$
then 

$$
\exists C \in (0,\infty), \ \sup_n || \sqrt{n} \
r(\hat{\theta_n}, \theta)|| B( \overline{\phi} ) = C(f(\cdot,\cdot)) < \infty. \eqno(6.10)
$$
 As a corollary: for all values $ p  = \const \in [1,\infty) $

$$
 \sup_n | \sqrt{n} \
r(\hat{\theta_n}, \theta)|_p \le C_1 \ p/\overline{\phi}^{-1}(p) < \infty, \eqno(6.11)
$$
where  we used the ordinary notation: for arbitrary r.v. $ \zeta $

$$
|\zeta|_p \stackrel{def}{=} \left[ {\bf E} |\zeta|^p  \right]^{1/p}.
$$

\vspace{3mm}

{\bf 7. Some examples.} \\

\vspace{3mm}

{\bf Example 7.1.} Spherical unimodal distributions. \\

\vspace{3mm}

We consider the following first example (and other examples) in order
to illustrate the precision  of  the theorems 4.1, 4.2 and 5.1. \par

 Let $ q = \const \ge 2, \ X = R^m, \ \mu $ be an usually Lebesgue measure,
$ x \in R^m \ \Rightarrow \ |x| = (x,x)^{1/2}; \
 R(y), y \in [0, \infty) \  $ be twice continuous differentiable
strictly  positive:

$$
\inf_{y \in [0,1]} R(y) > 0, \ \inf_{y \in [1, \infty)} y^q R(y) > 0,
$$
{\it slowly varying} as $ y \to \infty $ functions  such that the function
 $ y \to \ y^q R (y), \ y \ge 0 $ is strictly monotonically increasing. \par
 Let us introduce the following density function

$$
f_0(x) = C(q,m,R) \ \exp \left( - |x|^q \ R(|x|) \right),
$$
where $ C(q,m,R) $ is a norming  constant:

$$
\int_{R^m}  f_0(|x|) \ dx  = 1.
$$

 We take as a parametric set $ \Theta = X = R^m; $ choose $ \theta_0 := 0, $
and define the {\it family} of a densities of a view (shift family):

$$
f(x, \theta) = f_0(|x \ - \ \theta|), \ \theta \in \Theta = R^m. \eqno(7.1)
$$

 Recall that the observation (observations) $ \xi $ has (have) the density of
distribution $ f_0(|x|). $ \par

 It follows from the {\it unimodality} of the density function that the
MLE of the parameter $ \theta $ coincides with the observation  $ \xi: $ \par

$$
\hat{\theta} = \xi.  \eqno(7.2)
$$

\vspace{2mm}

{\bf A. Upper bound.} \par

\vspace{2mm}

 It follows after some computations on the basis of the theorem 5.1 that
(using the classical results from the theory of slowly, or regular
varying functions  functions)  (see \cite{Seneta1}, pp. 41 \ - \ 53)  that
for the function $ \phi(\lambda) = |\lambda|^q \ R(|\lambda|) $ the
Young \ - \ Fenchel transform has a following asymptotic: as $  \lambda \to
 \ \infty \ \Rightarrow $

$$
\left( |\lambda|^q \ R(|\lambda|) \right)^* \sim C \ |\lambda|^p /
R \left(|\lambda|^{p \ - \ 1} \right),
$$
where as usually $ 1/p + 1/q = 1. $ As long as $ q \ge 2, $ we conclude
that $ p \in (1, 2]. $ \par
 We obtain on the basis of theorem 5.1: $ v \ge 1 \  \Rightarrow $

$$
W(v) \le \exp \left[ - C_3 \ v^p \ /R \left( v^{p \ - \ 1} \right] \right).
\eqno(7.3)
$$

\vspace{2mm}

{\bf B. Low bound.} \par

\vspace{2mm}

We get using the explicit representation (7.2) and passing to the polar coordinates:

$$
W(v) = {\bf P} (|\xi| > v) = C(q,m,R) \int_{x: |x| > v} \
\exp \left( - |x|^q \ R(|x|) \right) \ dx =
$$

$$
C_9(q,m,R) \int_v^{\infty} y^{m \ - \ 1} \ \exp \left( - y^q \ R(y) \right)
\ dy \ge
$$

$$
C_{10} \ \exp \left( - v^q \ R(v) / C_{11} \right). \eqno(7.4)
$$

 Notice that the upper (7.3) and low bounds (7.4) exponential coincides if for instance $ p = q = 2 $ and $ R = \const  \ $ (the  "Gaussian case").\par
Analogously may be considered a more general case of the classical MLE estimations. \par

\vspace{3mm}

{\bf Example 7.2.} "Smooth"  sample. \par

\vspace{3mm}

 We suppose here that all the conditions of the theorem 6.1. are satisfied.\par
 It follows from the formula (6.8) that

$$
W_n(v) \le \exp \left( - n \ \phi^*(C_1 \ v/\sqrt{n} ) \right). \eqno(7.5)
$$
Assume that the variable $ v $ belongs to the following "zone": for some nonrandom positive constant $ C < \infty $

$$
v \le C_2 \ \sqrt{n} \eqno(7.6)
$$
(a "big" zone of great deviations). Substituting into (7.5) and taking into account
the behavior of the function  $ \phi^* = \phi^*(\lambda) $
 we obtain in the considered zone the estimation:

$$
W_n(v) \le  \exp \left( - C_3 \ v^2  \right). \eqno(7.7)
$$

 On the other hand, we observe that from the CLT for MLE estimations that for
 each {\it fixed} positive value $ v:$
$$
\overline{W}(v) \ge \lim_{n \to \infty} W_n(v) \ge \exp
\left( - C \ v^2 \right).
$$

\vspace{3mm}

{\bf Example 7.3.} Heavy tails of distributions. \par

\vspace{3mm}

We consider here the sample of a volume $ n $ from  the standard 
one \ - \ dimensional Cauchy distribution: $ X = R^1, \ \theta \in R^1, \theta_0 = 0, $

$$
f(x,\theta) = \frac{\pi^{-1}}{1 + (x \ - \ \theta)^2 }.
$$
 It is easy to calculate that

$$
\phi(\lambda) \asymp C_1 \ \lambda^2, \  |\lambda| \le C_2;
$$

$$
\phi(\lambda) \asymp C_3 \ |\lambda|, \  |\lambda| \ge C_2.
$$
 More fine considerations as in the theorem 6.1 based on the exponential
and power bounds for random fields maximum distribution based on the monograph
\cite{Ostrovsky2}, chapter 3, see also \cite{Ostrovsky3} show us that

$$
W_n(v) \le \exp \left(  - C_4 \ v^2\right), \ v \le C_5;
$$

$$
W_n(v) \le  C_6 /v, \ v \ge C_5.
$$

Therefore,

$$
\sup_n W_n(v) \le C_7/v, \ v \ge C_8. \eqno(7.8)
$$

On the other hand,

$$
\sup_n W_n(v) \ge W_1(v) = {\bf P} ( |\xi(1)| > v) \ge
C_9/v, \ v \ge C_{10}, \eqno(7.9)
$$
which coincides with upper bound (7.8) up to multiplicative constant.\par

 At the same result is true for the symmetric stable distributions with the shift parameter $ \theta: \ f(x, \theta) =  f(x - \theta), \ \theta \in R^1. $
 In detail,  let $ \{ f(\cdot, \cdot) \} $ be again the one \ - \ dimensional
shift family of densities with characteristical  functions

$$
\int_{-\infty}^{\infty} e^{i t x} \ f(x, \theta) \ dx =
e^{i t \theta  \ - \ |t|^{\alpha}}, \ \alpha \in (1,2).
$$
 Using at the same arguments we obtain the following bilateral inequality:

$$
C_1(\alpha)/v^{\alpha} \le  \overline{W}(v) \le C_2(\alpha)/v^{\alpha}, \ v \ge 1.
$$

\vspace{3mm}

{\bf Example 7.4.} Scale parameter.\\

\vspace{3mm}

 Let here $ \{ \xi(i) \}, \ i = 1,2,\ldots,n $ be a sample from  the one \ - \
dimensional distribution $ N(0, \theta), \ \theta > 0, X = R^1,
\theta_0 = 1. $\par

The theorem 6.1 gives us the following estimation:

$$
\overline{W}(v) \le \exp \left( - C \ v \right), v > 1. \eqno(7.10)
$$

The MLE $ \hat{\theta}_n $ has an explicit view:

$$
\hat{\theta}_n = n^{-1} \sum_{i=1}^n [\xi(i)]^2.
$$

 The distribution of $ \hat{\theta} $ coincides,  up to multiplicative constant,
with the known $ \chi^2 $ distribution with $ n $ degree of freedom. \par
We can see by means of this consideration that

$$
\overline{W}(v) \ge \exp \left( - C \ v \right), v > 1; \eqno(7.11)
$$
and moreover for all values $ n $

$$
W_n(v) \ge C_1(n) \exp  \left( - C_2(n) \ v  \right), \ v \ge C_3(n).    \eqno(7.12)
$$

 At the same result is true for exponential distribution, indeed, when

$$
f(x,\theta) = \theta^{-1} \ \exp \left( - x/\theta \right);
$$

$ X = R^1_+, \ \theta > 0, \ \theta_0 = 1. $ \par

Notice that in this case the value $ \lambda_0 $ from the definition (2.6) 
is  finite.\par

\vspace{4mm}

{\bf Remark 7.1} \\

\vspace{3mm}

 Note that the case of the so \ - \ called  {\it penalized} modification of the
MLE estimation (PMLE) may be considered analogously. See for definition
and first results in the nonasymptotic risk estimations in the PMLE
( \cite{Spokoiny1} ) and reference therein.\par

\end{document}